\documentclass[10pt,a4paper]{amsart}
\date{August 04, 2016}      
\usepackage{hyperref}
\usepackage{latexsym,amsmath,amsfonts,amscd,amssymb,graphics}
\usepackage[all]{xy}
\usepackage[utf8]{inputenc}

\theoremstyle{plain}  
\newtheorem{theorem}{Theorem}[section]
\newtheorem*{theorem*}{Theorem}

\newtheorem{corollary}[theorem]{Corollary}

\newtheorem{proposition}[theorem]{Proposition}

\theoremstyle{definition}
\newtheorem{definition}[theorem]{Definition}

\theoremstyle{remark}

\newtheorem{remark}[theorem]{Remark}

\newtheorem*{claim*}{Claim}

\numberwithin{equation}{section}

\renewcommand{\leq}{\leqslant}

\renewcommand{\geq}{\geqslant}

\renewcommand{\setminus}{\smallsetminus}

\newcommand{\R}{\mathbb{R}}

\newcommand{\Z}{\mathbb{Z}}
\newcommand{\C}{\mathbb{C}}

\newcommand{\PP}{\mathbb{P}}

\newcommand{\cM}{\mathcal{M}}
\newcommand{\cN}{\mathcal{N}}

\newcommand{\cR}{\mathcal{R}}
\newcommand{\cS}{\mathcal{S}}

\newcommand{\SO}{\mathrm{SO}}
\newcommand{\Or}{\mathrm{O}}
\newcommand{\Sp}{\mathrm{Sp}}

\newcommand{\PSp}{\mathrm{PSp}}

\newcommand{\U}{\mathrm{U}}

\DeclareMathOperator{\Jac}{Jac}

\DeclareMathOperator{\tr}{tr}
\DeclareMathOperator{\rk}{rk}

\DeclareMathOperator{\Hom}{Hom}

\DeclareMathOperator{\Id}{Id}

\DeclareMathOperator{\codim}{codim}

\newcommand{\Aut}{\operatorname{Aut}}

\let\oldmarginpar\marginpar
\renewcommand\marginpar[1]{\oldmarginpar{\tiny\bf\begin{flushleft} #1
\end{flushleft}}}

\begin{document}

\title[Quadric bundles and their moduli spaces]{Quadric bundles and their moduli spaces}

\author[A. G. Oliveira]{Andr\'e G. Oliveira}

\thanks{
Author partially supported by CMUP (UID/MAT/00144/2013), by the Projects EXCL/MAT-GEO/0222/2012 and PTDC/MAT-GEO/2823/2014 and also by the Post-Doctoral fellowship SFRH/BPD/100996/2014. These are funded by FCT (Portugal) with national (MEC) and European structural funds (FEDER), under the partnership agreement PT2020. Support from U.S. National Science Foundation grants DMS 1107452, 1107263, 1107367 ``RNMS: GEometric structures And Representation varieties'' (the GEAR Network) is also acknowledged. This note is partially based on the talk given at the Workshop NS@50: Fifty Years of the Narasimhan-Seshadri Theorem, celebrating the 50 years of the publication of \cite{narasimhan-seshadri:1965}, which was held at the Chennai Mathematical Institute (India), between 5 and 9 of October of 2015.}



\maketitle

\begin{flushright}
{\small{\it To M. S. Narasimhan and C. S. Seshadri,\\
on the occasion of the 50th anniversary of their Theorem.}} \\ \vspace*{0.25cm}
\end{flushright}

\vspace{.4cm}

\begin{abstract}
Quadric bundles on a compact Riemann surface $X$ generalise orthogonal bundles and arise naturally in the study of the moduli space of representations of $\pi_1(X)$ in $\Sp(2n,\R)$. We prove some basic results on the moduli spaces of quadric bundles over $X$ of arbitrary rank and then survey deeper results in the rank $2$ case.
\end{abstract}

\section{Introduction}

One can argue that one of the starting points on the study of moduli spaces of vector bundles over Riemann surfaces was achieved by the Narasimhan and Seshadri Theorem \cite{narasimhan-seshadri:1965}, which relates degree zero, stable vector bundles on a compact Riemann surface $X$ with unitary representations of $\pi_1(X)$. Indeed, as usually happens when groundbreaking results appear, this Theorem opened the doors to a all new world in Mathematics, giving origin to an explosion of research on moduli spaces of vector and principal bundles over projective curves. 
One of the opened doors was that of decorated vector bundles, that is, tuples consisting of vector bundles together with sections of associated bundles. One example is provided by \emph{quadric bundles} over compact Riemann surfaces, whose moduli spaces this note is devoted to. In the same spirit as the Narasimhan--Seshadri Theorem, these objects occur naturally (via Higgs bundles) by considering representations of the fundamental group of a compact Riemann surface not in $\U(n)$ but in the real symplectic group $\Sp(2n,\R)$. In addition, quadric bundles are natural generalisations of orthogonal vector bundles, being also a first case to consider if one wants to study decorated bundles whose sections do not live in a linear space, but rather in a higher tensor bundle. After describing these motivations to study quadric bundles, we prove some basic results on their moduli spaces for arbitrary rank. Then the remaining part of the paper, corresponding to Sections \ref{sec:4} and \ref{sec:5}, is a brief survey on some deeper results obtained in \cite{gothen-oliveira:2012,oliveira:2015} in the case of rank $2$.

\section{Quadric bundles}

\subsection{Definition as generalisation of orthogonal bundles}

Fix a compact Riemann surface $X$ of genus $g\geq 2$ and an arbitrary holomorphic line bundle $L\to X$.
An $L$-twisted orthogonal vector bundle over $X$ is a pair $(V,q)$ consisting of a vector bundle $V$ and an $L$-valued nowhere degenerate, symmetric, quadratic form $q:V\otimes V\to L$. Quadric bundles are natural generalisations of orthogonal vector bundles obtained by removing the non-degeneracy condition.

\begin{definition}
An \emph{$L$-quadric bundle} over $X$ is a pair $(V,\gamma)$ consisting of a holomorphic vector bundle $V$ on $X$ and a non-zero holomorphic section $\gamma$ of $S^2V^*\otimes L$.
\end{definition}

We shall say \emph{quadric bundle} instead of $L$-quadric bundle, whenever $L$ is our fixed line bundle and not a specific one. As usual, the \emph{type} of a quadric bundle $(V,\gamma)$ is $(n,d)$ if the rank and degree of $V$ are $n$ and $d$ respectively. We consider the section $\gamma$ as a map $\gamma:V\to V^*\otimes L$ which is symmetric, meaning $\gamma^t=\gamma\otimes\Id_L$, where $\gamma^t$ denotes the transpose or dual map.
No restrictions on the rank of $\gamma$ are imposed. If it is an isomorphism $V\cong V^*\otimes L$ then it is an $L$-valued nowhere degenerate symmetric quadratic form, thus $(V,\gamma)$ is an $L$-twisted orthogonal vector bundle.

Quadric bundles are sometimes also called \emph{conic bundles} or \emph{quadratic pairs} in the literature. In particular, this happens in the papers \cite{gothen-oliveira:2012,oliveira:2015} by the author where they were named quadratic pairs. But the term ``quadric bundles'', used in \cite{giudice-pustetto:2014}, is indeed more adequate since it is more specific. Moreover it emphasises the fact that these can be seen as bundles of quadrics, since for each $p\in X$ the map $\gamma$ restricted to the fibre $V_p$ defines a bilinear symmetric form, hence a quadric in $\mathbb{P}^{n-1}$. When $n=2$, the term conic bundle is then perfectly adequate also.

Quadric bundles are an example of a decorated vector bundle, meaning a bundle together with extra data, usually provided by sections of associated bundles. Other examples are Higgs bundles \cite{hitchin:1987,simpson:1992,bradlow-garcia-prada-gothen:2005}, Bradlow pairs \cite{bradlow:1991,thaddeus:1994} or holomorphic triples \cite{bradlow-garcia-prada-gothen:2004 triples}. In fact, one can argue that \emph{quadric bundles are a natural first case to consider if one wants to study decorated bundles whose associated sections are not linear but live instead in a higher tensor bundle}. 

Next we shall see that quadric bundles also arise naturally from the study of the moduli space of representations of $\pi_1(X)$ in the real symplectic group $\Sp(2n,\R)$, via Higgs bundles, in the same spirit as the Narasimhan--Seshadri Theorem.

\subsection{Quadric bundles and representations $\pi_1(X)\to\Sp(2n,\R)$}\label{subsec:NHC-Sp(n,R)}

After Ramanathan's generalisation  of the Narasimhan--Seshadri Theorem to compact reductive Lie groups in \cite{ramanathan:1975}, a natural question was to ask for a further generalisation to non-compact groups. That was achieved, for complex reductive Lie groups, after the introduction of \emph{Higgs bundles} by N. Hitchin in \cite{hitchin:1987}, through a combination of results of Hitchin \cite{hitchin:1987}, Donaldson \cite{donaldson:1987}, Corlette \cite{corlette:1988} and Simpson \cite{simpson:1992}.  More generally, given a real reductive, connected, Lie group $G$, there is the notion of $G$-Higgs bundle, which reduces to that of a $G^\C$-principal bundle when $G$ is compact. We will not define it, since we will only need it for $\Sp(2n,\R)$, defined below. We refer for example to \cite{bradlow-garcia-prada-gothen:2005} for the details. The following fundamental theorem --- the \emph{non-abelian Hodge correspondence} --- is then the extension of the Narasimhan--Seshadri Theorem as well as Ramanathan's Theorem to any such group $G$. We state it for $G$ semisimple, noting that there is a similar correspondence for the case $G$ reductive, by replacing $\pi_1(X)$ by the universal central extension.

Let $\cR_c(G)=\Hom^{\mathrm{red}}(\pi_1(X),G)/G$ be the space of equivalence classes of reductive representations of $\pi_1(X)$ in $G$, modulo the conjugation by $G$, and whose topological class of the underlying (flat) bundle is $c\in\pi_1(G)$. Let also $\cM_c(G)$ be the moduli space of polystable $G$-Higgs bundles on $X$, of the same topological type $c\in\pi_1(G)$. The non-abelian Hodge correspondence states the following:

\begin{theorem}[\cite{hitchin:1987,donaldson:1987,corlette:1988,simpson:1992,garcia-prada-gothen-mundet:2008}]\label{fundamental correspondence for semisimple G}
If $G$ is semisimple, then there is a natural correspondence between $\cM_c(G)$ and $\cR_c(G)$, which induces a homeomorphism $\cM_c(G)\cong\cR_c(G)$, for any $c\in\pi_1(G)$.
This correspondence comes from the fact that a $G$-Higgs bundle over $X$ is polystable if and only if it corresponds to a reductive representation of $\pi_1(X)$ in $G$. 
\end{theorem}

Let $K$ denote the canonical line bundle of $X$. Let $\Sp(2n,\R)$ be the real symplectic group i.e. the group of automorphisms of $\R^{2n}$ preserving a symplectic form. The definition of $G$-Higgs bundle in this case reads as follows. 

\begin{definition}
An \emph{$\Sp(2n,\R)$-Higgs bundle} is a triple $(V,\beta,\gamma)$ where $V$ is a holomorphic rank $n$ vector bundle, $\beta$ a section of $S^2V\otimes K$ and $\gamma$ a section of $S^2V^*\otimes K$.
\end{definition}

Hence one can think of $\beta$ as a holomorphic map $\beta:V^*\to V\otimes K$ such that $\beta^t\otimes \mathrm{Id}_K=\beta$ and likewise for $\gamma:V\to V^*\otimes K$.
The topological type of an $\Sp(2n,\R)$-Higgs bundle is the degree $d$ of $V$ --- this is the \emph{Toledo invariant}. There is a \emph{Milnor-Wood type of inequality} stating that $\cM_d(\Sp(2n,\R))$ is empty unless
\begin{equation}\label{eq:MWSp}
|d|\leq n(g-1).
\end{equation}

\begin{remark}
The Toledo invariant and the Milnor-Wood inequality are not exclusive of the case of $\Sp(2n,\R)$. They occur precisely when the group is of hermitian type \cite{biquard-garcia-prada-rubio:2015,bradlow-garcia-prada-gothen:2005}.
\end{remark}

Suppose then that $|d|\leq n(g-1)$. 
If we are interested in studying the connected components of $\cM_d(\Sp(2n,\R))$, then a standard tool is the \emph{Hitchin function} 
\[f:\cM_d(\Sp(2n,\R))\to\R,\qquad f(V,\beta,\gamma)=\|\beta\|_{L^2}^2+\|\gamma\|_{L^2}^2=\int_X\tr(\beta\beta^{*,h})+\int_X\tr(\gamma^{*,h}\gamma).\]
Here $h$  is the metric on $V$ which provides the Hitchin-Kobayashi correspondence between polystable
$\Sp(2n,\R)$-Higgs bundles and solutions to the \emph{Hitchin's equations} \cite{hitchin:1987,garcia-prada-gothen-mundet:2008,garcia-prada-gothen-mundet:2013}. The Uhlenbeck compactness theorem implies that $f$ is proper, which in turn says that the study of the connected components of $\cM_d(\Sp(2n,\R))$ reduces to the study of the components of the subvarieties of local minima of $f$. This entire paragraph extends in fact to any real reductive group.

It turns out that if  $0<d<n(g-1)$, then $(V,\beta,\gamma)$ represents a local minimum of $f$ if and only if $\beta\equiv 0$ (see \cite{garcia-prada-gothen-mundet:2013}), hence if and only if it is a $K$-quadric bundle. If $n(1-g)<d<0$ then it is $\gamma$ that vanishes; in any case $\cM_d(\Sp(2n,\R))$ and $\cM_{-d}(\Sp(2n,\R))$ are isomorphic by the duality $(V,\beta,\gamma)\mapsto(V^*,\gamma,\beta)$. 
Hence we see here that $K$-quadric bundles appear naturally in the study of $\cM_d(\Sp(2n,\R))$, therefore also on the space $\cR_d(\Sp(2n,\R))$, for non-zero and non-maximal Toledo invariant.

Indeed the proof of connectedness of $\cM_d(\Sp(2n,\R))$ for $0<|d|<n(g-1)$, has been carried out, for $n=2$, in \cite{garcia-prada-mundet:2004} and \cite{gothen-oliveira:2012}. The techniques in these papers are different, but both use $K$-quadric bundles. The method of \cite{gothen-oliveira:2012} was also applied to study the case of $\SO_0(2,3)$-Higgs bundles ($\SO_0(2,3)$ denotes the identity component of $\SO(2,3)$) --- hence representations $\pi_1(X)\to\SO_0(2,3)$ --- through the analysis of $L_0K$-quadric bundles, where $L_0$ is any degree $1$ line bundle on $X$.

\section{Moduli of quadric bundles}

\subsection{Semistability and non-emptiness of the moduli}
In this section we prove some basic results on the moduli space of quadric bundles for arbitrary rank. These moduli spaces were first considered by T. G\'omez and I. Sols in \cite{gomez-sols:2000}, for rank up to three. Their existence for arbitrary rank follows from the general construction of A. Schmitt in \cite{schmitt:2008}.

The following is the $\alpha$-semistability condition. It is equivalent to the $\delta$-semistability condition deduced in \cite{giudice-pustetto:2014} (see also \cite{gomez-sols:2000}), with $\alpha$ and $\delta$ related by $\alpha=\frac{d-2\delta}{n}$. We choose to use the parameter $\alpha$ to be consistent with \cite{garcia-prada-gothen-mundet:2013,gothen-oliveira:2012,oliveira:2015}.

Given a subbundle $V'\subseteq V$, denote by $V'^\perp$ the kernel of the projection $V^*\to V'^*$.
\begin{definition}\label{def:semistability-rkn}
Let $(V,\gamma)$ be an $L$-quadric bundle of type $(n,d)$ and let $\alpha\in\R$. Then:
\begin{itemize}
\item $(V,\gamma)$ is \emph{$\alpha$-semistable} if $\alpha\leq d/n$ and the following conditions hold:
\begin{enumerate}
\item for any proper vector subbundle $V'\subsetneq V$, 
\begin{enumerate}
\item if $\gamma(V')\not\subseteq V'^\perp\otimes L$, then $\deg(V')\leq d+\alpha(\rk(V')-n)$.
\item if $\gamma(V')\subseteq V'^\perp\otimes L$, then $\deg(V')\leq\frac{d}{2}+\alpha\left(\rk(V')-\frac{n}{2}\right)$.
\item if $\gamma(V')=0$, then $\deg(V')\leq\alpha\rk(V')$.
\end{enumerate}
\item for any strict filtration $0\subsetneq V'\subsetneq V''\subsetneq V$ such that $\gamma(V')\subseteq V''^\perp\otimes L$ and $\gamma(V'')\not\subseteq V''^\perp\otimes L$, we have $\deg(V')+\deg(V'')\leq d+\alpha(\rk(V')+\rk(V'')-n)$.
\end{enumerate}
\item $(V,\gamma)$ is \emph{$\alpha$-stable} if it is $\alpha$-semistable and strict inequalities hold in (1) and (2) above.
\item $(V,\gamma)$ is \emph{$\alpha$-polystable} if it is $\alpha$-semistable and the following conditions hold:
\begin{enumerate}
\item for any proper vector subbundle $V'\subsetneq V$, 
\begin{enumerate}
\item if $\gamma(V')\not\subseteq V'^\perp\otimes L$ and $\deg(V')=d+\alpha(\rk(V')-n)$, then there is a vector subbundle $V''\subsetneq V$ such that $V\cong V'\oplus V''$ and $\gamma=\left(\begin{smallmatrix}\gamma' & 0 \\0 & 0\end{smallmatrix}\right)$ with $\gamma'\in H^0(X,S^2V'^*\otimes L)$ non-zero.
\item if $\gamma(V')\subseteq V'^\perp\otimes L$ and $\deg(V')=d/2+\alpha(\rk(V')-n/2)$, then there is a vector subbundle $V''\subsetneq V$ such that $V\cong V'\oplus V''$ and $\gamma=\left(\begin{smallmatrix}0 & \gamma' \\ -\gamma'^t & 0\end{smallmatrix}\right)$ with $\gamma'\in H^0(X,V'^*\otimes V''^*\otimes L)$, non-zero. 
\item if $\gamma(V')=0$ and $\deg(V')=\alpha\rk(V')$, then there is a vector subbundle $V''\subsetneq V$ such that $V\cong V'\oplus V''$ and $\gamma=\left(\begin{smallmatrix}0 & 0 \\0 & \gamma'\end{smallmatrix}\right)$ with $\gamma'\in H^0(X,S^2V''^*\otimes L)$ non-zero.
\end{enumerate}
\item if there is a strict filtration in the conditions of (2) of the definition of $\alpha$-semistability, such that $\deg(V')+\deg(V'')= d+\alpha(\rk(V')+\rk(V'')-n)$, then $V\cong V'\oplus V''/V'\oplus V/V''$ and $\gamma=\left(\begin{smallmatrix}0 & 0 & -\gamma'^t\\0 & \gamma'' & 0 \\ \gamma' & 0 & 0 \end{smallmatrix}\right)$ with $\gamma'\in H^0(X, V'^*\otimes V/V''\otimes L)$ and $\gamma''\in H^0(X,S^2(V''/V')^*\otimes L)$, not both zero.
\end{enumerate}
\end{itemize}
\end{definition}

In general a quadric bundle can be (semi,poly)stable for a value of the parameter and unstable for other value. However, there are certain quadric bundles which are strictly semistable or polystable for any $\alpha$. We now specify these cases.

\begin{definition}\label{def:alphaindependent}
A quadric bundle $(V,\gamma)$ of type $(n,d)$ is said to be \emph{$\alpha$-independently strictly semistable} if one of the following two cases occur:
\begin{itemize}
\item both $n,d$ are even and $V$ has a proper subbundle $V'\subsetneq V$, such that $\deg(V')=d/2,\ \rk(V')=n/2$ and such that $\gamma(V')\subseteq V'^\perp\otimes L$.
\item  $V$ admits a strict filtration $0\subsetneq V'\subsetneq V''\subsetneq V$ of proper subbundles, such that $\rk(V')+\rk(V'')=n$, $\deg(V')+\deg(V'')=d$ and moreover $\gamma(V')\subseteq V''^\perp\otimes L$ and $\gamma(V'')\not\subseteq V''^\perp\otimes L$.
\end{itemize}
Analogously, one defines an \emph{$\alpha$-independently strictly polystable} quadric bundle.
\end{definition}

The first item of this definition should be compared with (1)(b) of the of the definition of $\alpha$-semistability in Definition \ref{def:semistability-rkn}. The second item should be compared with (2) of the same definition.

\vspace{.5cm}

Set once and for all $d_L=\deg(L)$. Let $\cN_\alpha(n,d)$ be the moduli space of $\alpha$-polystable $L$-quadric bundles over $X$ of type $(n,d)$, as constructed in \cite{gomez-sols:2000,schmitt:2008}.

Denote by $\rk(\gamma)$ the rank of  $\gamma:V\to V^*\otimes L$.

\begin{proposition}\label{prop:bounds}
Fix $\alpha$. If $(V,\gamma)$ is an $\alpha$-semistable quadric bundle of type $(n,d)$, then \[n\alpha\leq d\leq\frac{\rk(\gamma)d_L}{2}+(n-\rk(\gamma))\alpha.\]
In particular, if $\alpha> d_L/2$, then $\cN_\alpha(n,d)$ is empty for any $(n,d)$. If $\alpha\leq d_L/2$, then $\cN_\alpha(n,d)$ is empty whenever $d>nd_L/2$;
\end{proposition}
\proof
The first inequality follows directly from the $\alpha$-semistability condition. If $\gamma$ is generically non-degenerate then $\det(\gamma)$ is a non-zero section of $\det(V)^{-2}L^n$, hence $d\leq\frac{n}{2}d_L$ setting the second inequality in this case. Suppose now that $1\leq\rk(\gamma)\leq n-1$. Let $N\subseteq V$ be the kernel of $\gamma$ and $I\subseteq V^*$ be the subbundle such that $I\otimes L$ is the vector bundle generated by the image of $\gamma$.
Then $\gamma$ induces a non-zero map $\tilde\gamma:V/N\to I\otimes L$ of maximal rank, so the existence of the non-zero section $\det(\tilde\gamma)$ of $\det((V/N))^{-1}\det(I) L^{\rk(\gamma)}$ implies
\begin{equation}\label{eq:Milnor-Wood1}
-d+\deg(N)+\deg(I)+\rk(\gamma)d_L\geq 0.
\end{equation}
By $\alpha$-semistability, we have $\deg(N)\leq(n-\rk(\gamma))\alpha$ and also $\deg(I)\leq(n-\rk(\gamma))\alpha-d$, since $\gamma(I^\perp)=0$ by the symmetry of $\gamma$. From \eqref{eq:Milnor-Wood1} it follows that $-2d+2(n-\rk(\gamma))\alpha+\rk(\gamma)d_L\geq 0$, proving the second inequality.
The last sentences are obvious.
 \endproof

\begin{remark}
This result is consistent with the general statement of Theorem 1.1 of \cite{biquard-garcia-prada-rubio:2015}.
\end{remark}

We assume henceforth that $\alpha\leq d_L/2$, and want $n$ and $d$ such that $\alpha\leq d/n\leq d_L/2$. For fixed $\alpha\leq d_L/2$ and $n$, the allowed values for the integer $d$ lie then between $n\alpha$ and $nd_L/2$. Dividing these values into $n$ equal parts, that is taking $n\alpha<\dots<rd_L/2+(n-r)\alpha<\dots<nd_L/2$, and if $rd_L/2+(n-r)\alpha\leq d\leq nd_L/2$ for some with $r=0,\ldots,n$, then the moduli space $\cN_\alpha(n,d)$ only contains quadric bundles whose rank of $\gamma$ is at least $r$. As we decrease the value of $d$, then the allowed ranks of $\gamma$ in $\alpha$-semistable quadric bundles of degree $d$ increase.

From now on we tacitly assume that $n\in\Z_+$, $d\in\Z$, $\alpha\in\R$ and $d_L\in\Z$ satisfy 
\begin{equation}\label{eq:bounds n d}
n\alpha\leq d\leq\frac{n}{2}d_L.
\end{equation}

\subsection{Critical values}

Suppose we vary $\alpha$ and obtain $\alpha'\leq d/n$. Then the moduli spaces $\cN_\alpha(n,d)$ and $\cN_{\alpha'}(n,d)$ remain isomorphic unless there exists some quadric bundle $(V,\gamma)$ which is $\alpha$-stable and $\alpha'$-unstable, or the other way around. A consequence is the existence of some value $\alpha_c$ between $\alpha$ and $\alpha'$ (possibly equal to one of them) such that $(V,\gamma)$ is \emph{strictly} $\alpha_c$-polystable. Such $\alpha_c$ is called a critical value. More precisely, a \emph{critical value} $\alpha_c\leq d/n$ is a value of the parameter for which equalities are possible in the inequalities of Definition \ref{def:semistability-rkn}. That is, $\alpha_c=d/n$ or is such that there is a natural number $n'<n$ and an integer $d'$ so that
\[\alpha_c=\frac{d-d'}{n-n'}\ \text{ or }\ \alpha_c=\frac{d-2d'}{n-2n'}\ \text{ or }\ \alpha_c=\frac{d'}{n'}\] or there are natural numbers $n'<n''<n$ and integers $d',d''$ so that
\[\alpha_c=\frac{d-d'-d''}{n-n'-n''}.\]
Clearly the critical values form a discrete set. We shall now see that they are indeed finitely many.

Given a quadric bundle which is strictly $\alpha$-polystable then either $\alpha$ is a critical value or the quadric bundle is $\alpha$-independently strictly polystable; cf. Definition \ref{def:alphaindependent}.

Set 
\[\alpha_m=d-\frac{n-1}{2}d_L\hspace{.5cm}\text{and}\hspace{.5cm}\alpha_M=\frac{d}{n}.\]

\begin{proposition}\label{prop:boundscritval}
If $\alpha_c$ is a critical value, then $\alpha_m\leq\alpha_c\leq\alpha_M$.
\end{proposition}
\proof
From Proposition \ref{prop:bounds} we know that if $(V,\gamma)\in\cN_\alpha(n,d)$ for some $\alpha<\alpha_m$, then $\gamma$ must have maximal rank, that is, it must be generically non-degenerate. Suppose there exists a critical value $\alpha_c<\alpha_m$, so that there must exist a strictly $\alpha_c$-polystable pair $(V,\gamma)$, which furthermore is not $\alpha$-independently polystable. Hence $(V,\gamma)$ must fall into one of the cases of $\alpha$-polystability of Definition \ref{def:semistability-rkn}. We shall see that there is no such $(V,\gamma)$. 
Cases (1)(a) and (1)(c) cannot occur because there $\rk(\gamma)<n$.
In case (1)(b), we have $V\cong V'\oplus V''$ and $\gamma=\left(\begin{smallmatrix}0 & \gamma' \\ -\gamma'^t & 0\end{smallmatrix}\right)$. Since both $V'$ and $V''$ are isotropic, then from the non-degeneracy of $\gamma$, we must have $\rk(V')=n/2$. So $(V,\gamma)$ is $\alpha$-independently strictly polystable, thus it cannot be of the required type. 
It remains to deal with case (2). Here $V\cong V'\oplus V''/V'\oplus V/V''$ and $\gamma=\left(\begin{smallmatrix}0 & 0 & -\gamma'^t\\0 & \gamma'' & 0 \\ \gamma' & 0 & 0 \end{smallmatrix}\right)$. Since $\gamma(V')\subset(V/V'')^*\otimes L$ and $\gamma(V/V'')\subset V'^*\otimes L$, the non-degeneracy of $\gamma$ forces $\rk(V')+\rk(V'')=n$ thus, again, $(V,\gamma)$ must be $\alpha$-independently strictly polystable, hence not of the required type.
\endproof

Since there are no critical values less than $\alpha_m$, the next result is obvious.

\begin{corollary}\label{cor:nevenisom}
If $\alpha_1,\alpha_2<\alpha_m$, then $\cN_{\alpha_1}(n,d)$ and $\cN_{\alpha_2}(n,d)$ are isomorphic. 
\end{corollary}

\subsection{Large parameter}

Now we look at the other extreme of the range of $\alpha$. The next propositions give a description of the moduli spaces $\cN_\alpha(n,d)$ for values of $\alpha$ slightly smaller than $\alpha_M=d/n$. Let $\alpha_M^-=\frac{d}{n}-\epsilon$, such that $\epsilon>0$ is small enough so that there are no critical values between $\alpha_M^-$ and $\alpha_M$. 


\begin{proposition}\label{prop:largeparamater}
Let $(V,\gamma)$ be a quadric bundle of type $(n,d)$. If it is $\alpha_M^-$-semistable then $V$ is semistable as a vector bundle.
\end{proposition}
\proof
Let $V'\subsetneq V$ be a vector subbundle of degree $d'$ and rank $n'$. Then either $\gamma(V')\not\subseteq V'^\perp\otimes L$ or $0\neq \gamma(V')\subseteq V'^\perp\otimes L$ or $\gamma(V')=0$.
In the first case, $\alpha_M^-$-semistability of $(V,\gamma)$ implies $d'\leq\frac{dn'}{n}+\epsilon(n-n')$, which yields $\mu(V')\leq\mu(V)$ by letting $\epsilon$ go to $0$. The other two cases are similar.
\endproof

Let $M(n,d)$ denote the moduli space of rank $n$ and degree $d$ vector bundles over $X$ and write $M^s(n,d)$ for the stable locus.

\begin{corollary}\label{cor:largemoduli}
There is a forgetful map $\pi:\cN_{\alpha_M^-}(n,d)\to M(n,d)$.
\begin{itemize}
\item If $d<\frac{n}{2}(d_L+1-g)$, then $\pi$ is surjective.
\item If $d<\frac{n}{2}(d_L+2-2g)$, then the restriction of $\pi$ to $\pi^{-1}(M^s(n,d))$ is a projective bundle, with fiber $\PP^N$ where $N=\left(-d+\frac{n}{2}(d_L+1-g)\right)(n+1)-1$.
\end{itemize}
\end{corollary}
\proof
From the preceding proposition, the morphism $\pi:\cN_{\alpha_M^-}(n,d)\to M(n,d)$, $(V,\gamma)\mapsto V$ is well-defined. Given a semistable vector bundle $V$, the fiber of $\pi$ over $V$ is the quotient $H^0(S^2V^*\otimes L)\setminus\{0\}/\Aut(V)$, where an automorphism $f$ of $V$ acts on $H^0(S^2V^*\otimes L)\setminus\{0\}$ as $\gamma\mapsto (f^t\otimes\Id_L)\gamma f$. 

The surjectivity condition follows from the fact that $S^2V^*\otimes L$ is semistable and from Riemann-Roch, recalling that $\rk(S^2V^*)=\frac{n(n+1)}{2}$ and $\deg(S^2V^*)=-d(n+1)$. 

If $V$ is stable then $\Aut(V)\cong \C^*$, so the fiber is the projective space $\PP H^0(S^2V^*\otimes L)$. Moreover, $H^0(S^2V^*\otimes L)$ has constant dimension, equal to $\left(-d+\frac{n}{2}(d_L+1-g)\right)(n+1)$, whenever $d<\frac{n}{2}(d_L+2-2g)$.
\endproof

Remark \ref{rmk:last parameter-moduli bundles} below states how $\pi$ behaves when $\frac{n}{2}(d_L+1-g)\leq d\leq\frac{n}{2}d_L$.

\subsection{Maximal degree}

Here we shall consider the moduli spaces $\cN_\alpha(n,d)$ for $d=\frac{n}{2}d_L$, the maximal allowed value, according to \eqref{eq:bounds n d}.
Then $\alpha_m=\alpha_M$, so Corollary \ref{cor:nevenisom} says that all moduli spaces $\cN_\alpha(n,nd_L/2)$ are isomorphic, for every $\alpha<\alpha_M$.

Let $M_{\Or(n,\C)}$ denote the moduli space of rank $n$ orthogonal vector bundles, where here we are \emph{not} fixing the topological type which, if $n\geq 3$, is labeled by the first and second Stiefel-Whitney classes.

\begin{proposition}\label{prop:max degree - orthogonal}
If $d_L$ is even, then $\cN_\alpha(n,nd_L/2)$ is isomorphic to $M_{\Or(n,\C)}$, for any $\alpha<\alpha_M$.
\end{proposition}
\proof
Let $(V,\gamma)\in\cN_\alpha(n,nd_L/2)$. Then $\gamma$ has maximal rank. But since now $\deg(V)=\deg(V^*\otimes L)$, then $\gamma:V\to V^*\otimes L$ must in fact be a symmetric isomorphism. Since $d_L$ is even, we can choose a square root $L'$ of $L$. Hence, if we define $W=V\otimes L'^{-1}$ and $q=\gamma\otimes \Id_{L'^{-1}}$, then $(W,q)$ is an orthogonal vector bundle. 
Let us see that $(V,\gamma)$ is $\alpha$-semistable as a quadric bundle if and only if $(W,q)$ is semistable as an orthogonal bundle. Recall that $(W,q)$ is semistable if every isotropic subbundle of $W$ has non-negative degree. Since $\alpha<\alpha_M=\alpha_m$ we can, by Corollary \ref{cor:nevenisom}, assume $\alpha=\alpha_M^-=d/n-\epsilon$ as in the previous subsection.

If $(V,\gamma)$ is $\alpha_M^-$-semistable, then by Proposition \ref{prop:largeparamater},  $W$ is semistable as a vector bundle, hence $(W,q)$ is semistable as an orthogonal vector bundle.

Conversely, if $(W,q)$ is semistable as an orthogonal vector bundle, then Proposition 4.2 of \cite{ramanan:1981} states that $W$ is also semistable, hence so is $V$. So for any subbundle $V'\subsetneq V$, with $\deg(V')=d'$ and $\rk(V')=n'$, 
\begin{equation}\label{eq:d'n'dn}
d'/n'\leq d/n.
\end{equation}
Note that $\gamma(V')\neq 0$ because $\ker(\gamma)=0$.
Suppose $\gamma(V')\not\subseteq V'^\perp\otimes L$. Then, we must have \begin{equation}\label{eq:ss1}
d'\leq d+\alpha_M^-(n'-n).
\end{equation} due to \eqref{eq:d'n'dn}.
Suppose $\gamma(V')\subseteq V'^\perp\otimes L$. Then, again \eqref{eq:d'n'dn} implies 
\begin{equation}\label{eq:ss2}
d'\leq \frac{d}{2}+\alpha_M^-\left(n'-\frac{n}{2}\right),
\end{equation} since $n'\leq n/2$.
Finally, take a strict filtration $0\subsetneq V'\subsetneq V''\subsetneq V$ such that $\gamma(V')\subseteq V''^\perp\otimes L$ and $\gamma(V'')\not\subseteq V''^\perp\otimes L$. Let $\deg(V')=d'$, $\deg(V'')=d''$, $\rk(V')=n'$ and $\rk(V'')=n''$. Then
\begin{equation}\label{eq:d'n'dn2}
d'/n'\leq d/n\hspace{0.5cm} \text{ and }\hspace{0.5cm} d''/n''\leq d/n
\end{equation} and, furthermore, $n'+n''\leq n$ from the non-degeneracy of $\gamma$. Therefore \eqref{eq:d'n'dn2} implies 
\begin{equation}\label{eq:ss3}
d'+d''\leq d+\alpha_M^-(n'+n''-n).
\end{equation}
We conclude from \eqref{eq:ss1}, \eqref{eq:ss2} and \eqref{eq:ss3} that $(V,\gamma)$ is $\alpha_M^-$-semistable, completing the proof.
\endproof

The existence of this isomorphism is an example of the \emph{Cayley correspondence}; cf. \cite{bradlow-garcia-prada-gothen:2005,biquard-garcia-prada-rubio:2015}. It shows that, if $n\geq 3$ and $d_L$ is even, then $\cN_\alpha(n,nd_L/2)$ has at least $2\times 2^{2g}$ connected components.

\begin{remark}\label{rmk:last parameter-moduli bundles}
From Proposition \ref{prop:max degree - orthogonal} and from \cite{serman:2008} it follows that the map $\pi:\cN_{\alpha_M^-}(n,nd_L/2)\to M(n,nd_L/2)$ defined in Corollary \ref{cor:largemoduli} is an embedding, for $d_L$ even. On the other hand, if $\frac{n}{2}(d_L+1-g)\leq d<\frac{n}{2}d_L$, the determination of the image of $\pi$ is a Brill-Noether problem.
\end{remark}

\begin{remark}
If $d_L$ is odd then $L$ does not admit a square root, hence Proposition \ref{prop:max degree - orthogonal} does not apply. But everything else works (of course $n$ must be even), so one can say that $\cN_\alpha(n,nd_L/2)$ is isomorphic to the moduli space of $L$-twisted orthogonal vector bundles of rank $n$ and degree $nd_L/2$, where again we do not fix the topological type.
\end{remark}

\section{Connectedness}\label{sec:4}

From now on we assume that $n=2$. Accordingly, we write $\cN_\alpha(d)$ instead of $\cN_\alpha(2,d)$, with $d$ such that $2\alpha\leq d\leq d_L$. The remaining part of the article is a survey of the results obtained in \cite{gothen-oliveira:2012} and in \cite{oliveira:2015}, where all the details can be found. No new results are presented.
In fact we will be more interested in the case of non-maximal degree, so we actually assume \[2\alpha\leq d<d_L.\]

\subsection{Wall crossing}

Note that the conditions numbered with (2) of $\alpha$-semistability and of $\alpha$-polystability in Definition \ref{def:semistability-rkn} do not occur in rank $2$. Notice also that the $\alpha$-independent semistability (see Definition \ref{def:alphaindependent}) occurs precisely when one has a line subbundle $F\subset V$ such that $\deg(F)=d/2$ and $\gamma(F)\subseteq F^\perp\otimes L$.
In this case it is easy to conclude that, with the potential exception of $d/2$, all the critical values are integers. More precisely, the critical values are equal to $d/2$ or to $\alpha=\alpha_k=\lfloor d/2\rfloor+k$, for some $0\leq k\leq d-\lfloor d/2\rfloor-\lfloor d_L/2\rfloor.$
 By definition, on each open interval between consecutive critical values, the $\alpha$-semistability condition does not vary, hence the corresponding moduli spaces are isomorphic.
If $\alpha_k^+$ denotes the value of any parameter between the critical values $\alpha_k$ and $\alpha_{k+1}$, we can write without ambiguity 
$\cN_{\alpha_k^+}(d)$
for the moduli space of $\alpha_k^+$-semistable quadric bundles for any $\alpha$ between $\alpha_k$ and $\alpha_{k+1}$. Likewise, define $\cN_{\alpha_k^-}(d)$, with $\alpha_k^-$ denoting any value between the critical values $\alpha_{k-1}$ and $\alpha_k$. With this notation we always have $\cN_{\alpha_k^+}(d)=\cN_{\alpha_{k+1}^-}(d)$.
Moreover,
\begin{equation}\label{N--S-=N+-S+}
  \cN_{\alpha_k^-}(d)\setminus\cS_{\alpha_k^-}(d)
  \cong\cN_{\alpha_k^+}(d)\setminus\cS_{\alpha_k^+}(d).
\end{equation}
where $\cS_{\alpha_k^+}(d)\subset\cN_{\alpha_k^+}(d)$ is the subvariety consisting of those quadric bundles which are
$\alpha_k^+$-semistable but $\alpha_k^-$-unstable and analogously for  $\cS_{\alpha_k^-}(d)\subset\cN_{\alpha_k^-}(d)$. The spaces $\cS_{\alpha_k^\pm}(d)$ encode the difference between the spaces $\cN_{\alpha_k^-}(d)$ and $\cN_{\alpha_k^+}(d)$ on opposite sides of the critical value $\alpha_k$. This difference is usually known as the \emph{wall-crossing} phenomena through $\alpha_k$. 

In Section 3 of \cite{gothen-oliveira:2012} we studied the spaces $\cS_{\alpha_k^\pm}(d)$ for any critical value. We proved, in particular, that if $d<d_L-g+1$, then $\dim \cN_\alpha(d)=3(d_L-d)+g-1$ and, for any $k$, 
\begin{equation}\label{eq:codimS>g-1}
\codim \cS_{\alpha_k^\pm}(d)>g-1.
\end{equation} 

\subsection{Connectedness}

In \cite[Theorem 5.1]{gothen-oliveira:2012}, we proved that $\cN_{\alpha_m^-}(d)$ is connected for any $d<d_L$. The proof is based on the study made in \cite{gothen-oliveira:2013} on the Hitchin fibration for twisted Higgs pairs, i.e. analogues of Higgs bundles but with any twisting line bundle, not necessarily the canonical line bundle $K$. In a little more detail, in \cite{gothen-oliveira:2013} we made a sufficiently explicit study of the Hitchin fibration $h$ for rank $2$ Higgs pairs with trivial determinant, enough to prove that every fiber of $h$ is connected. Of course this is known for the generic fiber, which is a Prym variety, and in fact the connectedness of every fiber follows from this whenever the twisting line bundle has degree at least $2g-2$ (so works for $K$). However, we proved the result for any twisting line bundle (of course of positive degree). This yielded the connectedness of $\cN_{\alpha_m^-}(d)$  --- see Section 4 of \cite{gothen-oliveira:2012}.

This result, together with \eqref{eq:codimS>g-1} and \eqref{N--S-=N+-S+}, is enough to prove the following; see \cite[Theorem 5.3]{gothen-oliveira:2012}.

\begin{theorem}[\cite{gothen-oliveira:2012}]\label{thm:connected moduli}
For any $d$ such that $d<d_L-g+1$ and any $\alpha\leq d/2$, the moduli spaces $\cN_\alpha(d)$ are connected and non-empty.
\end{theorem}

The proof of the connectedness of the moduli spaces of $\Sp(4,\R)$-Higgs bundles (hence also of reductive representations of $\pi_1(X)$ in $\Sp(4,\R)$; cf. Theorem \ref{fundamental correspondence for semisimple G}) with non-zero and non-maximal Toledo invariant, mentioned in subsection \ref{subsec:NHC-Sp(n,R)}, follows from Theorem \ref{thm:connected moduli} by taking $L=K$. If one takes $L=L_0K$ for some degree one line bundle $L_0$, we reach the same conclusion for $G=\SO_0(2,3)\cong\PSp(4,\R)$. In the later case, the topological type is determined not only by the Toledo invariant $d$, but requires also a second Stiefel-Whitney class $w\in H^2(X,\Z/2)\cong\Z/2$ of a certain orthogonal bundle. Notice that here the constrain $d<d_L-g+1$ can be ruled out:

\begin{theorem}[\cite{gothen-oliveira:2012}]\mbox{}
\begin{itemize}
\item For any $d$ such that $0<|d|<2g-2$ the moduli spaces $\cM_d(\Sp(4,\R))$ and $\cR_d(\Sp(4,\R))$ are connected.
\item For any $d$ such that $0<|d|<2g-2$ and any $w\in\Z/2$, the moduli spaces $\cM_{d,w}(\SO_0(2, 3))$ and $\cR_{d,w}(\SO_0(2,3))$ are connected.
\end{itemize}
\end{theorem}

As we have already mentioned, the first part of this theorem as been previously proved in \cite{garcia-prada-mundet:2004}.

\section{Geometry of the fixed determinant moduli space}\label{sec:5}

From now on we consider rank $2$ quadric bundles whose determinant of the underlying bundle is a fixed line bundle $\Lambda\to X$ of degree $d$. Denote by $\cN_\alpha(\Lambda)$ this moduli space. The details of the results presented here can be found in \cite{oliveira:2015}.

\subsection{Irreducibility and simply-connectedness of the smooth locus}

The connectedness results of the preceding section started with the study of the lowest extreme of $\alpha$, that is the study of $\cN_{\alpha_m^-}(d)$.
If we start with the highest possible extreme, by taking $\alpha_M^-=d/2-\epsilon$, we have Corollary \ref{cor:largemoduli} (which immediately applies to the fixed determinant situation) telling us what happens. 
So let $M(2,\Lambda)$ be the moduli space of rank $2$  vector bundles with fixed determinant $\Lambda\in\Jac^d(X)$. Then we have the forgetful map
\begin{equation}\label{eq:forgetmap}
\pi:\cN_{\alpha_M^-}(\Lambda)\to M(2,\Lambda),\qquad \pi(V,\gamma)=V.
\end{equation} 
Again denote by $M^s(2,\Lambda)$ the corresponding stable locus. So we have as in Corollary \ref{cor:largemoduli},
\begin{proposition}\label{prop:bundle-forget}
The restriction of $\pi$ in \eqref{eq:forgetmap} to $\pi^{-1}(M^s(2,\Lambda))$ is a projective bundle, with fiber $\PP^{3(d_L-d+1-g)-1}$, whenever $d<d_L+2-2g$.
\end{proposition}

It is important to note that the spaces $\cN_\alpha(\Lambda)$ are \emph{not smooth}. Indeed they are singular precisely at the points represented by strictly $\alpha$-polystable quadric bundles (due to $\alpha$-independent polystability these also may occur for $\alpha$ non-critical) and by those $\alpha$-stable quadric bundles which have other automorphisms than $\pm\Id_V$, called \emph{non-simple} quadric bundles.
 In \cite{oliveira:2015} we explicitly described the singular locus. In particular we concluded that $\cN_\alpha(\Lambda)$ is smooth out of a subspace of high codimension (precisely, $3g-3$ if $\alpha$ is non-critical and $g-1$ if $\alpha$ is critical) and that the smooth locus $\cN_\alpha^{sm}(\Lambda)$ is dense for any $\alpha$.

Furthermore, we described the space $\pi^{-1}(M^s(2,\Lambda))$, showing that it is contained in $\cN_\alpha^{sm}(\Lambda)$ and that the complement has codimension $g-1$. Using this description, Proposition \ref{prop:bundle-forget} and the fact that $M^s(2,\Lambda)$ is irreducible, we concluded that $\cN_{\alpha_M^-}(\Lambda)$ is irreducible as well. Then, using \eqref{eq:codimS>g-1} and \eqref{N--S-=N+-S+}, we proved the following:

\begin{proposition}
If $g\geq 2$ and $d<d_L+2-2g$ then $\cN_\alpha(\Lambda)$ is irreducible.
\end{proposition}

The irreducibility of $\cN_\alpha(\Lambda)$ has been proved before in \cite[Theorem 3.5]{gomez-sols:2000}, under the same conditions, using different techniques. 

Roughly, this is the basic (standard) technique we use in \cite{oliveira:2015}: start with the analysis of the projective bundle \eqref{eq:forgetmap} when $d<d_L+2-2g$ and then, with the information we have from $M^s(2,\Lambda)$, draw conclusions on the space $\cN^{sm}_{\alpha_M^-}(\Lambda)$. After, from \eqref{eq:codimS>g-1} and \eqref{N--S-=N+-S+}, we transfer our conclusions to $\cN^{sm}_\alpha(\Lambda)$ for any $\alpha$. Finally we might be able to conclude something about the entire moduli $\cN_\alpha(\Lambda)$.

This was again performed in the next result, using that $M^s(2,\Lambda)$ is simply-connected whenever $g\geq 3$ or $g=2$ and $d=\deg(\Lambda)$ is odd.

\begin{proposition}
Let $g\geq 3$ or $g=2$ and $d$ odd. If $d<d_L+2-2g$ then $\cN^{sm}_\alpha(\Lambda)$ is simply-connected for every $\alpha<\alpha_M$.
\end{proposition}

Observe that we cannot take the final step and conclude that $\cN_\alpha(\Lambda)$ is simply-connected, since compactifying a smooth space by a singular one may seriously affect the fundamental group. In any case, it would be interesting to determine $\pi_1(\cN_\alpha(\Lambda))$.

\subsection{Cohomology groups and a Torelli theorem}

Now we state our result on the cohomology groups of $\cN_\alpha^{sm}(\Lambda)$. First we recall some facts on Hodge structures.
A \emph{(pure) Hodge structure of weight $k$} is a pair $(V_\Z,V_\C)$ where $V_\Z$ is a free abelian group and $V_\C=V_\Z\otimes\C$ is its complexification, such that there is a decomposition into subspaces $V_\C=\bigoplus_{p+q=k}V^{p,q}$ with $\overline{V^{q,p}}\cong V^{p,q}$.
The Hodge structure $(V_\Z,V_\C)$ is \emph{polarised} if it admits a bilinear non-degenerate map $\theta:V_\Z\otimes V_\Z\to\Z$, alternating if $k$ is odd and symmetric otherwise, such that $V^{p,q}$ and $V^{p',q'}$ are orthogonal under $\theta_\C:V_\C\otimes V_\C\to\C$, unless $p'=q$ and $q'=p$, and moreover $i^{p-q}(-1)^{k(k-1)/2}\theta_\C(\alpha,\overline\alpha)>0$ for every $\alpha\in V^{p,q}$.

The cohomology of an $n$-dimensional smooth projective variety $Z$ (hence compact K\"ahler, with K\"ahler form $\omega$) has a pure Hodge structure of weight $k$ given by $H^k(Z,\C)=\bigoplus_{p+q=k}H^{p,q}(Z)$. This Hodge structure is polarised by the wedge product:
\begin{equation}\label{eq:polarisation of smooth projective variety}
\theta(\alpha,\beta)=\int_Z\alpha\wedge\beta\wedge\omega^{n-k}, \qquad \alpha,\beta\in H^k(Z,\Z).
\end{equation}
However, this may not be true anymore for varieties which are not smooth or projective. In these cases Deligne showed that we have to consider the generalised concept of mixed Hodge structure \cite{deligne:1970,deligne:1971,deligne:1974}. This is the case of the moduli spaces $\cN^{sm}_\alpha(\Lambda)$ which are smooth, quasi-projective, but not projective. Nevertheless, it can be seen, using Corollary 6.1.2 of \cite{arapura-sastry:2000}, that at least the first three cohomology groups of $\cN^{sm}_\alpha(\Lambda)$ have a \emph{pure} Hodge structure. So we will not really use in practice the notion of mixed Hodge structure.

From Lemma 6.1.1 and Corollary 6.1.2 of \cite{arapura-sastry:2000}, from Proposition \ref{prop:bundle-forget} and using the basic technique explained in the previous subsection, we proved the following (using ideas introduced by V. Mu\~noz in \cite{munoz:2009}):

\begin{proposition}[\cite{oliveira:2015}]\label{Prop:cohomology groups}
Let $g\geq 4$ and $d<d_L+2-2g$. Then, for every $\alpha<\alpha_M$, the torsion-free part of the first three cohomology groups of $\cN_\alpha^{sm}(\Lambda)$ is given by 
\[H^1(\cN^{sm}_\alpha(\Lambda),\Z)=0,\quad H^2(\cN^{sm}_\alpha(\Lambda),\Z)\cong\Z\oplus\Z,\quad H^3(\cN^{sm}_\alpha(\Lambda),\Z)\cong H^1(X,\Z).\]
Moreover, the Picard group of $\cN_\alpha^{sm}(\Lambda)$ is isomorphic to $\Z\oplus\Z$.
\end{proposition}

Finally, note that $H^1(X,\Z)$ is polarised by $\theta$ given in \eqref{eq:polarisation of smooth projective variety}. Moreover, one can naturally polarise $H^3(\cN^{sm}_\alpha(\Lambda),\Z)$ in such a way that, if $g\geq 5$, the isomorphism $H^3(\cN^{sm}_\alpha(\Lambda),\Z)\cong H^1(X,\Z)$ preserves the polarisations; cf. \cite[Proposition 5.8]{oliveira:2015}. Since the classical Torelli Theorem says that the pair $(H^1(X,\Z),\theta)$ determines the (isomorphism class of the) curve $X$, we conclude that a Torelli type theorem also holds for $\cN^{sm}_\alpha(\Lambda)$ and hence for $\cN_\alpha(\Lambda)$:

\begin{theorem}[\cite{oliveira:2015}]
Let $X$ and $X'$ be smooth projective curves of genus $g,g'\geq 5$, $\Lambda$ and $\Lambda'$ line bundles of degree $d$ and $d'$ on $X$ and $X'$, respectively, and $L$ and $L'$ line bundles of degree $d_L$ and $d_{L'}$ on $X$ and $X'$, respectively. Assume $d<d_L+2-2g$ and $d'<d_{L'}+2-2g$. Let $\cN_{X,\alpha}(\Lambda)$  be the moduli space of $\alpha$-polystable $L$-quadric bundles on $X$, and define similarly $\cN_{X',\alpha}(\Lambda')$.
If $\cN_{X,\alpha}(\Lambda)$ and $\cN_{X',\alpha}(\Lambda')$ are isomorphic as abstract projective varieties, then so are $X$ and $X'$. 
\end{theorem}

\begin{remark}
From this theorem one can in fact prove that the same holds for the non-fixed determinant moduli spaces $\cN_\alpha(d)$; see \cite[Section 5.4]{oliveira:2015}.
\end{remark}

It would be very interesting to generalise to higher ranks the results presented in Sections \ref{sec:4} and \ref{sec:5}. One of the main obstacles is the description of the spaces $\cS_{\alpha_k^\pm}(d)$  involved in the wall-crossing through critical points, even just enough to estimated their codimension as in \eqref{eq:codimS>g-1}. This description becomes much more challenging for higher rank, also because of the complexity of the stability condition.
\vspace{.5cm}

\normalsize

\noindent
      \textbf{Andr\'e Oliveira} \\
      Centro de Matem\'atica da Universidade do Porto, CMUP\\
      Faculdade de Ci\^encias, Universidade do Porto\\
      Rua do Campo Alegre 687, 4169-007 Porto, Portugal\\ 
      email: andre.oliveira@fc.up.pt

\vspace{.2cm}
\noindent
\textit{On leave from:}\\
 Departamento de Matem\'atica, Universidade de Tr\'as-os-Montes e Alto Douro, UTAD \\
Quinta dos Prados, 5000-911 Vila Real, Portugal\\ 
email: agoliv@utad.pt

\end{document}